\documentclass[journal]{IEEEtran}
\usepackage{cite}
\usepackage{amssymb}
\usepackage{graphicx}
\usepackage{amsfonts}
\usepackage{lineno}
\usepackage{amsmath}
\usepackage{array}
\usepackage{subeqnarray}
\usepackage{cases}

\def\^{\hat}
\def\~{\tilde}

\def\3h{{3\over 2}}

\def\v#1{{\bf#1}}

\def\E{{\v E}}

\def\W{{\v W}}

\def\H{{\v H}}

\def\F{{\v F}}

\def\N{{\v N}}

\def\ep{\epsilon}

\def\p{\partial}

\def\div{\nabla\cdot}

\def\curl{\nabla\times}
\def\grad{\nabla}

\def\eqn#1$${\eqno{{\rm #1}}$$}

\def\ep{\epsilon}

\def\~{\tilde}

\def\^{\hat}

\def\bm#1{\mbox{{\boldmath${#1}$}}}

\def\d#1{\overline{\overline{#1}}}

\def\XXint#1#2#3{{\setbox0=\hbox{$#1{#2#3}{\int}$}
     \vcenter{\hbox{$#2#3$}}\kern-.5\wd0}}

\allowdisplaybreaks
\begin{document}
\title{Three Numerical Eigensolvers for 3-D Cavity Resonators Filled With Anisotropic and Nonconductive Media}
\author{\IEEEauthorblockN{Wei Jiang and Jie Liu}
\thanks{
This work was supported by the National Natural Science Foundation of China under Grant 61901131, the Natural Science Foundation of Guizhou Minzu University under Grant GZMU[2019]YB07 and the China Postdoctoral Science Foundation under Grant 2019M662244 \emph{(Corresponding author: Wei Jiang).}}
\thanks{W. Jiang is with the School of Mechatronics Engineering, Guizhou Minzu University, Guiyang 550025, China (e-mail: jwmathphy@163.com).
}
\thanks{J. Liu is with the Postdoctoral Mobile Station of Information and Communication
Engineering, School of Informatics, and Institute of Electromagnetics
and Acoustics, Xiamen University, Xiamen 361005, China (e-mail: liujie190484@163.com).
}}


\maketitle
\begin{abstract}
This paper mainly investigates the classic resonant cavity problem with anisotropic and nonconductive media, which is a linear vector Maxwell's eigenvalue problem. The finite element method based on edge element of the lowest-order and standard linear element is used to solve this type of 3-D closed cavity problem. In order to eliminate spurious zero modes
in the numerical simulation, the divergence-free condition supported by Gauss' law is enforced in a weak sense. After the finite element discretization, the generalized eigenvalue problem with a linear constraint condition needs to be solved. Penalty method, augmented method and projection method are applied to solve this difficult problem in numerical linear algebra. The advantages and disadvantages of these three computational methods are also given in this paper. Furthermore, we prove that the augmented method is free of spurious modes as long as the anisotropic material is not magnetic lossy. The projection method based on singular value decomposition technique can be used to solve the resonant cavity problem. Moreover, the projection method {cannot} introduce any spurious modes. At last, several numerical experiments are carried out to verify our theoretical results.

\end{abstract}


\maketitle
\IEEEdisplaynontitleabstractindextext
\IEEEpeerreviewmaketitle

\section{Introduction}
\IEEEPARstart{M}{icrowave} resonant cavity is an important passive device in microwave engineering. It has many applications
in many projects, such as particle accelerator, microwave oven and microwave filter. The microwave resonant cavity problem usually needs to solve the eigenmodes of  source-free Maxwell's equations. If the cavity is filled with inhomogeneous media and/or has a complex geometry, then finding its resonant modes by an analytical method is impossible. In order to get the resonant mode, numerical methods must be applied to solve the problem. The main numerical methods include finite element method, finite difference method and boundary element method.

Solving 3-D closed cavity problem will introduce many spurious modes if the numerical method
{cannot} preserve the physical property of electromagnetic field. These spurious modes are usually divided into two kinds. One is spurious nonzero mode and
the other one is spurious zero mode. In general, the spurious zero mode in numerical results is caused by the neglection of divergence-free condition, and the introduction of the spurious nonzero mode in numerical results is the result of improper discretization method.
It is our goal to design a numerical method that can eliminate these two kinds of spurious modes together. We now know edge element method can remove all spurious nonzero modes in solving electromagnetic resonant cavity problems. There are many references on this subject, such as Lee and Mittra\cite{Lee1992}, Wang and Ida \cite{Wang1992}, Pichon and Razek \cite{Pichon1992}, and so on. However, the edge element method {cannot} eliminate spurious zero modes because it {cannot} guarantee the solenoidal properties of electric and magnetic flux densities in the whole cavity.

The main difficulty of solving resonant cavity problems is how to enforce the divergence-free condition given by Gauss's law in electromagnetics. For the first time, F. Kikuchi \cite{Kikuchi} introduces a Lagrange multiplier to  enforce the divergence-free condition in a weak sense, and propose a mixed finite element method (MFEM) to solve 3-D empty cavity problem. As a consequence, Kikuchi's method is free of all the spurious modes, including spurious zero modes. Based on Kikuchi's idea, Jiang \emph{et al.} \cite{Jiang2016,Jiang12019} make use of MFEM to solve 3-D resonant cavity problems with anisotropic media, and successfully remove spurious zero and nonzero modes together. However, the MFEM supported in \cite{Jiang2016,Jiang12019} cannot deal with 3-D closed cavity problems filled with the anisotropic media, which is both electric and magnetic lossy. On the basis of the reference \cite{Jiang2016}, Liu \emph{et al.} \cite{Liu2017} give a two-grid vector discrete scheme for 3-D cavity problems with lossless media. The scheme given in \cite{Liu2017} is free of all spurious modes and is very efficient if we just need to know the first few physical modes. Using edge element and linear element, Jiang \emph{et al.} \cite{Jiang2017} successfully solve 3-D closed cavity problem filled with fully conductive media in the whole cavity. The numerical method given in \cite{Jiang2017} can also remove the spurious zero and nonzero modes together.

This paper continues to study the microwave cavity problem filled with anisotropic and nonconductive media. After eliminating the electric field in source-free Maxwell's equations, one can get a vector Maxwell's eigenvalue problem for magnetic field. This problem can be transformed into a corresponding variational formulation by using Green's formulaes. The edge basis functions of the lowest-order and standard nodal basis functions of linear element are used to discretize the variational formulation. Finally, we need to solve the generalized eigenvalue problem with a linear constraint condition, which is a very difficult problem in numerical linear algebra. In order to overcome this difficult problem, penalty method, augmented method and projection method reduce it to the generalized eigenvalue problem without any constraint. In addition, the advantages and disadvantages among these three computational methods are also given in the paper.

The outline of the paper is as follows. The governing equations and finite element discretization of 3-D resonant cavity problem are given in Section \uppercase\expandafter{\romannumeral2}. In Section \uppercase\expandafter{\romannumeral3}, we provide the penalty method, augmented method and projection method to solve the constrained generalized eigenvalue problem and discuss the advantages and disadvantages among these three numerical computational methods. In Section \uppercase\expandafter{\romannumeral4},
three numerical experiments are carried out to verify our theoretical results.

\section{Finite Element Discretization of 3-D Resonant Cavity Problem}
\subsection{Governing Equations for 3-D Resonant Cavity Problem}
Suppose that $\Omega$ is a bounded domain in $\mathbb{R}^{3}$, $\p\Omega$ is the boundary of $\Omega$ and $\^n$ is the outward normal unit vector on $\p\Omega$. Let $\ep_{0}$ and $\mu_{0}$ be the permeability and permittivity in vacuum, respectively. The relative permeability and permittivity tensor of an anisotropic medium are denoted by $\d{\mu}_{r}$ and $\d{\ep}_{r}$, respectively. The angular frequency of electromagnetic wave is denoted by $\omega$.

The governing equations of 3-D closed cavity problem are
the source-free Maxwell's equations of the first-order. After eliminating the electric field $\E$ in the source-free Maxwell's equations, one can get a second-order partial differential equations (PDEs) for the magnetic field $\H$:
\begin{subequations} \label{eq:3}
\begin{numcases}{}
  \curl\Big({\d{\ep}_{r}^{-1}}\curl\H\Big) =\Lambda{\d{\mu}}_{r}\H~~\text{in}~\Omega\label{eq:3a}\\
  \div\big({\d{\mu}}_{r}\H\big)= 0 ~~\text{in}~\Omega\label{eq:3b}\\
  \^n\times({\d{\ep}_{r}^{-1}}\curl\H)={\bf{0}}~~\text{on}~\p\Omega\label{eq:3c}\\
  \^n\cdot(\d{\mu}_{r}\H) = {0}~~\text{on}~\p\Omega\label{eq:3d},
\end{numcases}
\end{subequations}
where $\Lambda=\omega^2\ep_{0}\mu_{0}$ is the square of the wavenumber in vacuum. We would like to seek $(\Lambda,\H)$ with $\H\neq{\bf{0}}$ such that PDEs (\ref{eq:3}) holds. In electromagnetics, $(\Lambda,\H)$ with $\H\neq{\bf{0}}$ is called a physical resonant eigenmode in the resonant cavity. In mathematics, $(\Lambda,\H)$ with $\H\neq{\bf{0}}$ is called an eigenpair of PDEs (\ref{eq:3}), and $\Lambda$ and $\H$ are called the eigenvalue and eigenfunction of PDEs (\ref{eq:3}), respectively. In addition,
we know PDEs (\ref{eq:3}) only has a discrete point spectrum. Note that there may be several zero eigenmodes in PDEs (\ref{eq:3}), for details, please see \cite{Jiang2019p}.

If the anisotropic material is lossless, then $\d{\mu}_{r}$ and $\d{\ep}_{r}$ are both Hermitian \cite{Chew1990book}, that is
$\d{\mu}_{r}^{\dag}=\d{\mu}_{r}$ and $\d{\ep}_{r}^{\dag}=\d{\ep}_{r}$, where the superscript ${\dag}$ denotes the conjugate transposition. Moreover, when the anisotropic material is lossless, assuming that $\d{\ep}_{r}$ and $\d{\mu}_{r}$ are Hermitian positive definite since the anisotropic and lossless material in nature usually has this property. For the sake of simplicity, a Hermitian positive definite matrix $\bm{M}$ is denoted by $\bm{M}^{\dag}=\bm{M}>0$. In terms of the lossy characteristics of the anisotropic and nonconductive material, it is usually divided into the following four categories:
\begin{enumerate}
  \item Case 1: $\d{\ep}_{r}^{\dag}=\d{\ep}_{r}>0$ and $\d{\mu}_{r}^{\dag}=\d{\mu}_{r}>0$. The medium is lossless.
    \item Case 2: $\d{\ep}_{r}^{\dag}\neq\d{\ep}_{r}$ and $\d{\mu}_{r}^{\dag}=\d{\mu}_{r}>0$. The medium is electric lossy, but is not magnetic lossy.
  \item Case 3: $\d{\ep}_{r}^{\dag}=\d{\ep}_{r}>0$ and $\d{\mu}_{r}^{\dag}\neq\d{\mu}_{r}$. The medium is magnetic lossy, but is not electric lossy.
  \item Case 4: $\d{\ep}_{r}^{\dag}\neq\d{\ep}_{r}$ and $\d{\mu}_{r}^{\dag}\neq\d{\mu}_{r}$. The medium is both electric and magnetic lossy.
\end{enumerate}

Under Case 1, 2 and 3, the MFEM can deal with these types of the resonant cavity problems very well,
and it is free of all the spurious modes in numerical results \cite{Jiang2016,Jiang12019}.
However, the MFEM is not suitable for 3-D closed cavity problem under Case 4, because it is difficult to propose an appropriate mixed variational formulation. For the 3-D closed cavity problem under Case 4, the projection method introduced in this paper can deal with this problem very well, and the projection method can remove all the spurious modes.

\subsection{Finite Element Discretization}
It is well-known that the finite element method is a variational method, and only operates on the weak form of PDE, instead of the strong form of PDE. Hence, the corresponding weak form associated with PDEs (\ref{eq:3}) is given in the subsection. To get this weak form, we introduce the following Hilbert spaces over $\mathbb{C}$:
\begin{gather*}
     L^{2}(\Omega)=\big\{f: \int_{\Omega}|f(x,y,z)|^2 \mathrm {d}x\mathrm{d}y\mathrm{d}z<+\infty\big\}\\
    H^{1}(\Omega)=\big\{f\in{L^2(\Omega):\grad f\in{(L^2(\Omega))^3}}\big\}\\
    \H(\mbox{curl},\Omega)=\big\{\F\in{(L^2(\Omega))^3}: \curl{\F}\in{(L^2(\Omega))^3}\big\}.
\end{gather*}
Define the continuous sesquilinear forms:
\begin{eqnarray*}
\mathcal{A}:&&\H(\mbox{curl},\Omega)\times\H(\mbox{curl},\Omega)\rightarrow{\mathbb{C}}\\
&&(\H,\F)\rightarrow\int_{\Omega}{\d{\ep}_{r}^{-1}\curl\H\cdot\curl\F^{*}}d\Omega \\
\mathcal{M}:&& (L^2(\Omega))^3\times (L^2(\Omega))^3\rightarrow{\mathbb{C}}\\
&& (\H,\F)\rightarrow\int_{\Omega}{\d{\mu}_{r}\H\cdot\F^{*}}d\Omega\\
\mathcal{C}:&&\H(\mbox{curl},\Omega)\times{H^1(\Omega)}\rightarrow{\mathbb{C}}\\
&&(\H,q)\rightarrow\int_{\Omega}{\d{\mu}_{r}\H\cdot{\grad{q^{*}}}}d\Omega
\end{eqnarray*}
where the symbol $*$ stands for the complex conjugation of a given complex-valued function.

Using the Green's formulas, the weak form of PDEs (\ref{eq:3}) reads as:
\begin{subequations} \label{eqnt:4}
\begin{numcases}{}
\textrm{Seek~} \Lambda\in{\mathbb{C}},~\H\in{\H(\mbox{curl},\Omega)},~{\H}\neq\bf{0} \textrm{~such that}\nonumber\\
\mathcal{A}(\H,\F)= \Lambda \mathcal{M}(\H,\F),~\forall~\F\in{\H(\text{curl},\Omega)}\label{eqnt:4a}\\
 \mathcal{C}(\H,q) = 0,~\forall~q\in{H^1(\Omega)}\label{eqnt:4b}
\end{numcases}
\end{subequations}
Under Case 1, the eigenvalues $\Lambda$ are made up of countable nonnegative real numbers. Under Case 2, 3 and 4, the eigenvalues $\Lambda$ are made up of countable complex numbers. The physical interpretation is such a physical fact that there is no electromagnetic energy loss in the resonant cavity if the material is lossless and there exists electromagnetic energy loss in the resonant cavity provided that the material has a dielectric loss.

We now consider the conforming finite element discretization of (\ref{eqnt:4}). Let $\mathcal{T}_{h}$ be a regular tetrahedral mesh of the cavity $\Omega$. Here $h$ is the length of the longest edge in tetrahedral mesh $\mathcal{T}_{h}$. As usual, the edge element space $\W^{h}$ of the
lowest-order under the mesh $\mathcal{T}_{h}$ is used to approximate
the Hilbert space $\H(\mbox{curl},\Omega)$ and standard linear element space $S^{h}$ under the mesh $\mathcal{T}_{h}$ is used to approximate the Hilbert space $H^{1}(\Omega)$. From \cite{hiptmair2002}, we know $S^{h}\subsetneq H^{1}(\Omega)$ and $\W^{h}\subsetneq\H(\mbox{curl},\Omega)$.

The linear element space $S^{h}$ can be written as:
\begin{equation*}
    S^{h}=\big\{\phi:\phi|_{K}\in\textrm{span}\{L_{1}^{K},L_{2}^{K},L_{3}^{K},L_{4}^{K}\}\big\}
\end{equation*}
where $L_{i}^{K}~(i=1,2,3,4)$ are four local nodal basis functions on the tetrahedral element $K$ and of the form $a_{i}+b_{i}x+c_{i}y+d_{i}z$, where $a_{i},b_{i},c_{i},d_{i}$ are four constants. These four local basis functions are defined on the four vertices of the tetrahedral element $K$.

\begin{figure}[!t]
\centering
\includegraphics[width=4.5cm]{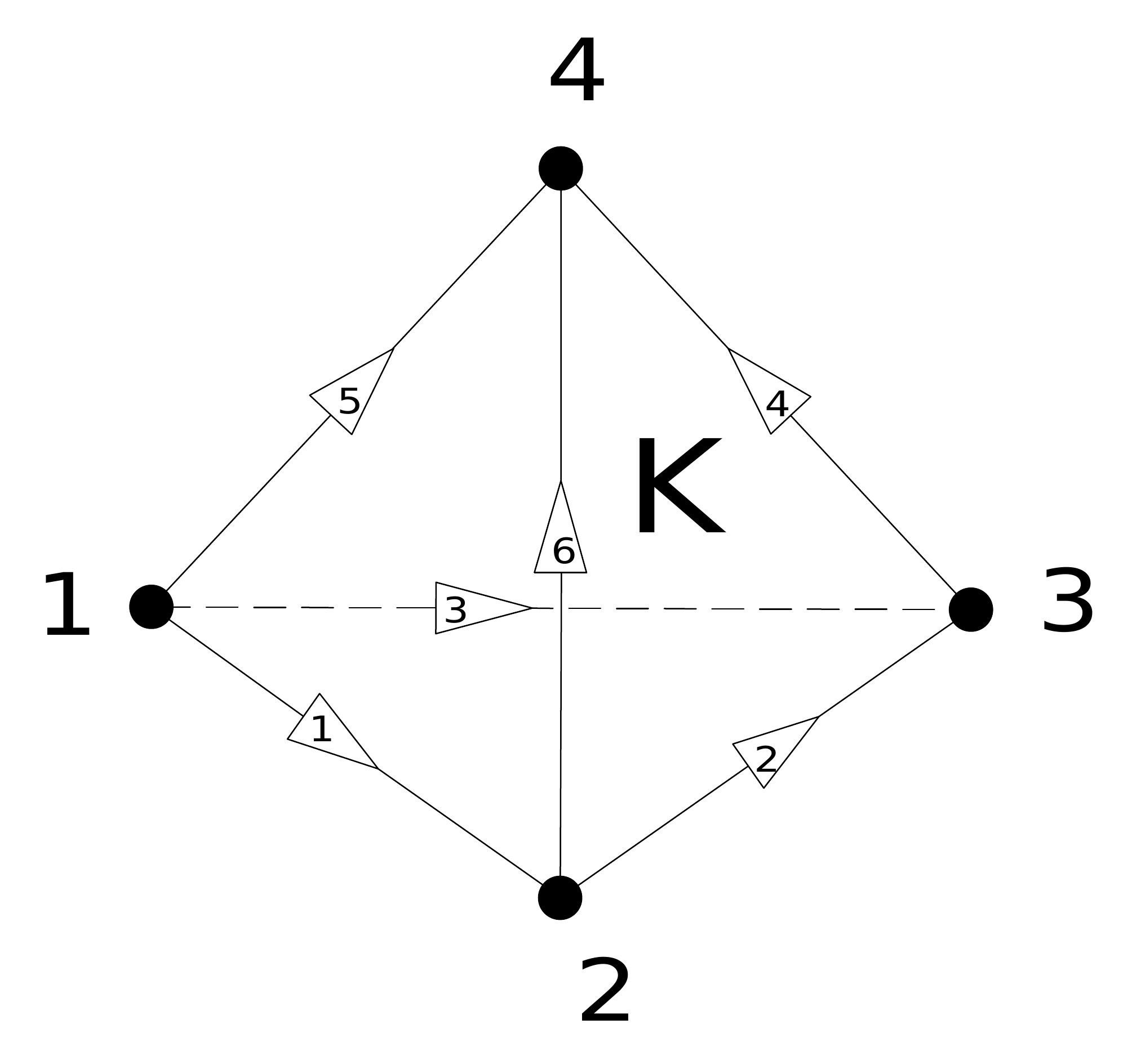}
\caption{Local nodal numbering for the element $K$ and the local reference direction for the edge are chosen by means of local nodal numbering.}\label{fig1}
\end{figure}

The edge element space $\W^{h}$ of the lowest-order can be written as:
\begin{equation*}
   \W^{h}=\big\{\F:\F|_{K}\in\textrm{span}\{\bm{N}_{1}^{K},\bm{N}_{2}^{K},\cdots,\bm{N}_{6}^{K}\}\big\}
\end{equation*}
where $\bm{N}_{i}^{K}~(i=1,2,\cdots,6)$ are six local edge basis functions on the tetrahedral element $K$ and $\bm{N}_{i}^{K}$ is of the form $\vec{\alpha}_{i}+\vec{\beta}_{i}\times\vec{r}$, where $\vec{\alpha}_{i}$ and $\vec{\beta}_{i}$ are two constant vectors and $\vec{r}$ is the position vector. The concrete expressions of $\bm{N}_{i}^{K}~(i=1,2,\cdots,6)$ are as follows:
\begin{gather*}
    \bm{N}_{1}^{K}=L_{1}^{K}\grad L_{2}^{K}-L_{2}^{K}\grad L_{1}^{K},~ \bm{N}_{2}^{K}=L_{2}^{K}\grad L_{3}^{K}-L_{3}^{K}\grad L_{2}^{K}\\
        \bm{N}_{3}^{K}=L_{1}^{K}\grad L_{3}^{K}-L_{3}^{K}\grad L_{1}^{K},~ \bm{N}_{4}^{K}=L_{3}^{K}\grad L_{4}^{K}-L_{4}^{K}\grad L_{3}^{K}\\
            \bm{N}_{5}^{K}=L_{1}^{K}\grad L_{4}^{K}-L_{4}^{K}\grad L_{1}^{K},~ \bm{N}_{6}^{K}=L_{2}^{K}\grad L_{4}^{K}-L_{4}^{K}\grad L_{2}^{K}
\end{gather*}
These six local vector basis functions are defined on the six edges of the tetrahedral element $K$. In Fig. \ref{fig1}, we give a local nodal numbering in the tetrahedral element $K$,
and specify the local reference direction for each edge in $K$.

Restricting (\ref{eqnt:4}) on $\W^{h}\times S^{h}$, we get
the discrete variational formulation associated with (\ref{eqnt:4}):
\begin{subequations} \label{eqnd2}
\begin{numcases}{}
\textrm{Seek~} \Lambda_{h}\in{\mathbb{C}},~\H_{h}\in{\W^{h}},~\H_{h}\neq{\bf{0}} \textrm{~such that}\nonumber\\
\mathcal{A}(\H_{h},\F)= \Lambda_{h} \mathcal{M}(\H_{h},\F),~\forall~\F\in{\W^{h}}\label{eqnd2a}\\
 \mathcal{C}(\H_{h},q) = 0,~\forall~q\in{S^{h}}\label{eqnd2b}
\end{numcases}
\end{subequations}
Here $\Lambda_{h}$ and $\H_{h}$ are an approximation of the exact eigenvalue $\Lambda$ and  exact eigenfunction $\H$ in (\ref{eqnt:4}), respectively.

Suppose that $S^{h}=\textrm{span}\big\{L_{i}\big\}_{i=1}^{m}$, where $L_{i}$ is the $i$-th global nodal basis function associated with the node $i$ and the integer $m$ {is} the number of the total nodes in $\mathcal{T}_{h}$. Assuming that
$\W^{h}=\textrm{span}\big\{\bm{N}_{i}\big\}_{i=1}^{n}$,
where $\bm{N}_{i}$ is the $i$-th global edge basis function associated
with the edge $i$ and the integer $n$ is the number of the total edges in $\mathcal{T}_{h}$.
Since $\H_{h}\in{\W^{h}}$, then it can be written as
\begin{equation}\label{Hspan}
    \H_{h}=\sum_{i=1}^{n}\xi_{i}\bm{N}_{i}.
\end{equation}

Finally, the discrete variational formulation (\ref{eqnd2}) can be reduced to the following generalized eigenvalue problem with a linear constraint:
\begin{subequations} \label{eigp2}
\begin{numcases}{}
\bm{A}\xi
=\Lambda_{h} \bm{M}\xi\label{eigp2a}\\
\bm{C}\xi={\bf{0}}\label{eigp2b}
\end{numcases}
\end{subequations}
where
\begin{gather*}
    \bm{A}=(a_{ik})\in{\mathbb{C}^{n\times n}},~\bm{M}=(m_{ik})\in{\mathbb{C}^{n\times n}},\\ ~\bm{C}=(c_{ik})\in{\mathbb{C}^{m\times n}},~\xi=[\xi_{1},~~\xi_{2},\cdots,\xi_{n}]^{T}\in{\mathbb{C}^{n}},\\
    a_{ik}=\mathcal{A}(\N_{k},\N_{i}),~m_{ik}=\mathcal{M}(\N_{k},\N_{i}),~
    c_{ik}=\mathcal{C}(\N_{k},L_{i}).
\end{gather*}
Solving (\ref{eigp2a}) directly and ignoring (\ref{eigp2b}) will introduce a lot of spurious zero modes. In order to eliminate these spurious zero modes, we must enforce the constraint (\ref{eigp2b}) in solving (\ref{eigp2a}). How to enforce the constraint (\ref{eigp2b}) in solving (\ref{eigp2a}) is a key problem. In next section, we shall deal with this troublesome problem. Once the eigenpair $(\Lambda_{h},\xi)$ is found from (\ref{eigp2}), then the numerical eigenvalue in (\ref{eqnd2}) is given by $\Lambda_{h}$ and the corresponding numerical eigenfunction $\H_{h}$ in (\ref{eqnd2}) is given by (\ref{Hspan}).

In the community of numerical linear algebra, the problem (\ref{eigp2}) is a constrained generalized eigenvalue problem. Obviously, its numerical computation is much more difficult than that of the generalized eigenvalue problem without any constraint. In Section III, the problem (\ref{eigp2}) is reduced to three types of generalized eigenvalue problem without any constraint. {It is} important to point out if there is no relation among these three matrices $\bm{A}$, $\bm{M}$ and $\bm{C}$,
then the constrained generalized eigenvalue problem (\ref{eigp2}) may have no solution.

\section{Three Numerical Solvers of the Constrained Generalized Eigenvalue Problem (\ref{eigp2})}
In this section, we first try to give a relation among the matrices $\bm{A}$, $\bm{M}$ and $\bm{C}$, and then support three numerical computational methods of solving (\ref{eigp2}). They are penalty method, augmented method and projection method, respectively.
\subsection{The Relation Among the Matrices $\bm{A}$, $\bm{M}$ and $\bm{C}$}

Let $\{A_{1},A_{2},\cdots,A_{m-1},A_{m}\}$ and $\{e_{1},e_{2},\cdots,e_{n-1},e_{n}\}$ be the sets consisting of all nodes and edges in $\mathcal{T}_{h}$ respectively, where $1,2,\cdots,m$ and $1,2,\cdots,n$ are the global labels of all nodes and edges in $\mathcal{T}_{h}$, respectively. Note that the global vector basis function $\bm{N}_{i}$ has a local direction associated with the edge $e_{i}$. If two vertices of the edge $e_{i}$ are $A_{i_{1}}$ and $A_{i_{2}}$,
then we state that the direction of $\bm{N}_{i}$ is from the node $A_{i_{1}}$ to the node $A_{i_{2}}$, where $i_{1}<i_{2}$.

In accordance with deRham-complex \cite{Bossavit1988}, $\grad{S^{h}}\subsetneq \W^{h}$ holds, where $\grad{S^{h}}=\{\grad f:\,\forall f\in{S^{h}}\}$. This implies that
\begin{equation}\label{imm}
    \grad{L_{i}}=\sum_{k=1}^{n}y_{ik}\N_{k},~~i=1,2,\cdots, m.
\end{equation}
The above formula (\ref{imm}) is also introduced in \cite{geus} and \cite{White2002}. It is easy to know $\dim(\grad{S^{h}})=m-1$ since $\sum_{i=1}^{m}\grad{L_{i}}={\bf{0}}$. Set
\begin{equation*}
   \bm{Y}= \begin{bmatrix}
      y_{11}&y_{12}&\cdots&y_{1n}\\
      y_{21}&y_{22}&\cdots&y_{2n}\\
      \vdots&\vdots&\vdots&\vdots\\
      y_{m1}&y_{m2}&\cdots&y_{mn}
    \end{bmatrix}=\begin{bmatrix}
    \bm{y}_{1}\\
    \bm{y}_{2}\\
    \vdots\\
    \bm{y}_{m}
    \end{bmatrix}=[\bm{d}_{1},\bm{d}_{2},\cdots,\bm{d}_{n}],
\end{equation*}
where $\bm{y}_{k}$ and $\bm{d}_{k}$ are $k$-th row and column vector in the matrix $\bm{Y}$, respectively.

In fact, the each entry $y_{ki}$ in $\bm{Y}$ is $-1$, $1$ or $0$, and $\bm{Y}$ is quite sparse.
Let us  consider the formula (\ref{imm}) under the case of an arbitrary tetrahedral element $K$ in $\mathcal{T}_{h}$ (see Fig. \ref{fig1}). It is easy to verify that the following formulas are valid:
\begin{gather*}
    \grad L_{1}^{K}=-\bm{N}_{1}^{K}-\bm{N}_{3}^{K}-\bm{N}_{5}^{K},~~\grad L_{2}^{K}=\bm{N}_{1}^{K}-\bm{N}_{2}^{K}-\bm{N}_{6}^{K}\\
    \grad L_{3}^{K}=\bm{N}_{2}^{K}+\bm{N}_{3}^{K}-\bm{N}_{4}^{K},~~\grad L_{4}^{K}=\bm{N}_{4}^{K}+\bm{N}_{5}^{K}+\bm{N}_{6}^{K}
\end{gather*}
Here we need to make use of the relation $\sum_{i=1}^{4}L_{i}^{K}=1$. Consider the $k$-th row $\bm{y}_{k}$ in the matrix $\bm{Y}$, which is related to the node $A_{k}$. If $A_{k}$ is not a vertex of the $i$-th edge $e_{i}$, then $y_{ki}=0$. Let us recall the basic concept of degree of a vertex in graph theory. The degree of a vertex is defined by the number of edges connecting it. The number of the nonzero entries of $\bm{y}_{k}$ is equal to the degree of $A_{k}$. Assuming that the degree of $A_{k}$ is $v$ and $\{e_{i_{1}},e_{i_{2}},\cdots,e_{i_{v}}\}$ have the common vertex $A_{k}$. If the direction of $e_{i_{s}}$ points to $A_{k}$,
then $y_{k,i_{s}}=1$, otherwise $y_{k,i_{s}}=-1$. Consider the $k$-th column vector $\bm{d}_{k}$ in the matrix $\bm{Y}$, which is related to the edge $e_{k}$. If $A_{i}$ is not a vertex of the edge $e_{k}$, then $y_{ik}=0$. Since each edge in a tetrahedral mesh only has two vertices, the column vector $\bm{d}_{k}$ only has two nonzero entries. Assuming that the initial and terminal points of the edge $e_{k}$ are $A_{i_{1}}$ and $A_{i_{2}}$ respectively, then $y_{i_{1}k}=-1$ and $y_{i_{2}k}=1$, where $i_{1}<i_{2}$. Obviously, the sparse matrix $\bm{Y}$ can be easily obtained by the mesh data in $\mathcal{T}_{h}$. The sparse matrix $\bm{Y}$ is usually called the directional connectivity matrix associated with a tetrahedral mesh $\mathcal{T}_{h}$.

The sum of all entries in the column vector $\bm{d}_{k}$ is equal to zero. This implies that the homogenous linear equation
\begin{equation}\label{holinear1}
    \bm{Y}^{\dag}\zeta={\bf{0}}
\end{equation}
has a special nonzero solution $\beta=[1,1,\cdots,1]^{T}\in{\mathbb{C}^{m}}$. An easy computation shows that $\textrm{rank}(\bm{Y})=\textrm{rank}(\bm{Y}^{\dag})=m-1$. Therefore all solutions of (\ref{holinear1}) form a linear space of one dimension, that is
\begin{equation}\label{linearsol}
    Null(\bm{Y}^{\dag})=\textrm{span}\{\beta\},
\end{equation}
where \emph{Null} stands for taking the nullspace of a matrix.

The matrix $\bm{Y}$ builds a relation among $\bm{A}$, $\bm{M}$ and $\bm{C}$. It can be proved that the following matrix identities are valid:
\begin{equation}\label{matrix12}
    \bm{Y}\bm{A}=\bm{O},~~\bm{C}=\bm{Y}\bm{M},
\end{equation}
where $\bm{O}$ is a null $m\times n$ matrix. In fact, for $\forall~1\leq i\leq m$ and $\forall~1\leq l\leq n$, we have
\begin{eqnarray*}
    (\bm{Y}\bm{A})_{il}&=&\sum_{k=1}^{n}y_{ik}a_{kl}=
    \sum_{k=1}^{n}y_{ik}\mathcal{A}(\N_{l},\N_{k})\\
    &=&\mathcal{A}(\N_{l},\sum_{k=1}^{n}y_{ik}\N_{k})=
    \mathcal{A}( \N_{l},\grad L_{i})=0,\\
   (\bm{Y}\bm{M})_{il}&=&\sum_{k=1}^{n}y_{ik}m_{kl}=
    \sum_{k=1}^{n}y_{ik}\mathcal{M}(\N_{l},\N_{k})\\
    &=&\mathcal{M}(\N_{l},\sum_{k=1}^{n}y_{ik}\N_{k})=
    \mathcal{M}( \N_{l},\grad L_{i})\\
    &=&\mathcal{C}( \N_{l},L_{i})=c_{il},
\end{eqnarray*}
where $(\bm{S})_{il}$ is the entry at $i$-th row and $l$-th column of the matrix $\bm{S}$.
{It is} important to emphasize that the matrix identities (\ref{matrix12}) are always valid for the medium under Case 1, 2, 3 and 4. Hence, we can obtain the matrix $\bm{C}$ by using the sparse matrices $\bm{Y}$ and $\bm{M}$, instead of the calculation by using the continuous sesquilinear forms $\mathcal{C}$ directly.

If the eigenvalue $\Lambda_{h}$ is nonzero in (\ref{eigp2a}),
then (\ref{eigp2b}) can be deduced from (\ref{eigp2a}).
As a matter of fact, $\bm{M}\xi=\Lambda_{h}^{-1}\bm{A}\xi$ can be derived from (\ref{eigp2a}), and then we get
$\bm{C}\xi=\bm{Y}\bm{M}\xi=\Lambda_{h}^{-1}\bm{Y}\bm{A}\xi={\bf{0}}$ by (\ref{matrix12}), which is just (\ref{eigp2b}). It is worthwhile to point out the number of zero eigenvalues in (\ref{eigp2a}) is equal to $\dim(\grad{S^{h}})$, and these zero eigenvalues are all spurious. The dimension of the linear space $\grad{S^{h}}$ is equal to $m-1$, which shows that the larger
the number of mesh nodes $m$ is, the more the number of spurious zero modes in (\ref{eigp2a}) is. Based on this reason, we need to remove these spurious zero modes by introducing an appropriate numerical method.

\subsection{Penalty Method}
Consider the following generalized eigenvalue problem:
\begin{equation}\label{numeig1}
    (\bm{A}+\alpha\bm{C}^{\dag}\bm{C})\xi=\Lambda_{h}' \bm{M}\xi,~~\|\xi\|_{2}=1
\end{equation}
where $\alpha$ is a penalty constant, which is required the user to set and $\|\varphi\|_{2}$ is the Euclidean norm of a given vector $\varphi$. The problem (\ref{numeig1}) is a generalized eigenvalue problem without any constraint, and it can be solved by the numerical software package ARPACK \cite{Lehoucq}.

The parameter $\alpha$ is usually set to a large positive real number.  The reason is as follows: Obviously, from (\ref{numeig1}), it is easy to deduce that
\begin{equation}\label{apprxomai1}
    \frac{1}{\alpha}\|\bm{A}\xi-\Lambda_{h} \bm{M}\xi\|_{2}=\|\bm{C}^{\dag}\bm{C}\xi\|_{2}.
\end{equation}
When $\alpha$ approaches positive infinity, we can obtain the following homogeneous linear equation from (\ref{apprxomai1}):
\begin{equation}\label{approa1}
    \bm{C}^{\dag}\bm{C}\xi={\bf{0}}.
\end{equation}
Multiplying both sides of (\ref{approa1}) on the left with ${\xi}^{\dag}$, then we can get $\|\bm{C}\xi\|_{2}=0$. As a result, the homogeneous linear equation (\ref{eigp2b}) is obtained again. Substituting (\ref{eigp2b}) into (\ref{numeig1}) gives $\bm{A}\xi=\Lambda_{h}'\bm{M}\xi$. However, if $\alpha$ takes sufficiently large, then this will lead to $\bm{A}+\alpha\bm{C}^{\dag}\bm{C}\approx\alpha\bm{C}^{\dag}\bm{C}$ because of the finite word length in a computer. The result is that the information of matrix $\bm{A}$ is completely submerged. Consequently, the numerical accuracy of the eigenvalues associated with the physical modes in (\ref{numeig1}) will become very poor \cite{Hayata}. Therefore, we {cannot} take a sufficiently large penalty parameter $\alpha$ in (\ref{numeig1}). It is worthwhile that choosing an appropriate penalty parameter $\alpha$ is important to the penalty method. Webb \cite{Webb1988} has studied this problem, and support a method to select this suitable parameter $\alpha$.

It can be seen that the eigenpair $(\Lambda_{h},\xi)$ of (\ref{eigp2}) is always the eigenpair $(\Lambda_{h},\xi)$ of (\ref{numeig1}). However, the eigenpair $(\Lambda_{h}',\xi)$ of (\ref{numeig1}) is not always the eigenpair of (\ref{eigp2}).
That is to say that penalty method will introduce many spurious modes in solving (\ref{numeig1}). Therefore, this method is not perfect. The penalty method is applicable to 3-D closed cavity problem under Case 1, 2, 3 and 4, except that it {cannot} remove all the spurious modes. If the eigenvalues in (\ref{numeig1}) are known, then how to choose  the eigenvalues with physical significance is an important problem. Here, we introduce two methods to identify these eigenvalues associated with physical modes:
\begin{enumerate}
  \item Assuming that $(\Lambda_{h}',\xi)$ with $\|\xi\|_{2}=1$ is an eigenpair of (\ref{numeig1}). The quantity $\|\bm{C}\xi\|_{2}$ can be used to identify the eigenvalues with physical significance. If $\|\bm{C}\xi\|_{2}$ is very small, then $\Lambda_{h}'$ is an eigenvalue corresponding to a physical mode. Otherwise, $\Lambda_{h}'$ is an eigenvalue associated with the spurious mode.
  \item The eigenvalues corresponding to physical eigenmodes will not change as the penalty constant $\alpha$ changes, but the eigenvalues corresponding to spurious modes will change as the penalty constant $\alpha$ changes. Therefore, under the same mesh, set $\alpha=\alpha_{1},\alpha_{2},\alpha_{3},\cdots$ and then solve (\ref{numeig1}) repeatedly, finally select the eigenvalues that remain unchanged in this processing. These unchanged eigenvalues are physical, whereas the changed eigenvalues are spurious.
\end{enumerate}

\subsection{Augmented Method}
Consider the generalized eigenvalue problem:
\begin{subequations} \label{numeig2}
\begin{numcases}{}
\bm{A}\xi+\bm{C}^{\dag}\zeta
=\Lambda_{h}' \bm{M}\xi\label{numeig2a}\\
\bm{C}\xi={\bf{0}}\label{numeig2b}
\end{numcases}
\end{subequations}
Obviously, the problem (\ref{numeig2}) can be rewritten as the following matrix form:
\begin{equation}\label{eigstt1}
    \left[
      \begin{array}{cc}
        \bm{A} & \bm{C}^{\dagger} \\
        \bm{C} & \bm{O} \\
      \end{array}
    \right]\left[
             \begin{array}{c}
               \xi \\
               \zeta \\
             \end{array}
           \right]
    =\Lambda_{h}'\left[
      \begin{array}{cc}
        \bm{M} & \bm{O} \\
        \bm{O} & \bm{O} \\
      \end{array}
    \right]\left[
             \begin{array}{c}
               \xi \\
               \zeta \\
             \end{array}
           \right].
\end{equation}
In the numerical linear algebra, the software package ARPACK \cite{Lehoucq} can be used to solve (\ref{eigstt1}).

It is clear that the necessary and sufficient condition for the equivalence of the eigenpair between (\ref{eigp2}) and  (\ref{eigstt1}) is
\begin{equation}\label{equiv1}
    \bm{C}^{\dag}\zeta={\bf{0}}.
\end{equation}
Assuming that $(\Lambda_{h}',[\xi;\zeta])$ is the eigenpair of (\ref{eigstt1}), where $[\xi;\zeta]$ is the corresponding eigenvector associated with $\Lambda_{h}'$ in (\ref{eigstt1}). Multiplying both sides of (\ref{numeig2a}) on the left with $\bm{Y}$, we obtain
\begin{equation}\label{auxi}
  (\bm{Y}\bm{A})\xi+(\bm{Y}\bm{C}^{\dag})\zeta
=\Lambda_{h}' (\bm{Y}\bm{M})\xi
\end{equation}
Substituting (\ref{matrix12}) and into (\ref{auxi}) gives
\begin{equation}\label{auxiw}
\bm{Y}\bm{C}^{\dag}\zeta=\Lambda_{h}' \bm{C}\xi={\bf{0}}.
\end{equation}
If (\ref{equiv1}) can be derived from (\ref{auxiw}),
then the eigenpair $(\Lambda_{h}',[\xi;\zeta])$ of (\ref{eigstt1}) is also an eigenpair $(\Lambda_{h}',\xi)$ of (\ref{eigp2}). Conversely, assuming that $(\Lambda_{h},\xi)$ is the eigenpair of (\ref{eigp2}).
If we take a vector $\zeta$ such that (\ref{equiv1}) is valid (Obviously, this is achievable), then the eigenpair $(\Lambda_{h},\xi)$ of (\ref{eigp2}) is an eigenpair $(\Lambda_{h},[\xi;\zeta])$ of (\ref{eigstt1}). This is to say that each eigenvalue of (\ref{eigp2}) is always the eigenvalue of (\ref{eigstt1}).

According to the above discussion, it can be concluded that the necessary and sufficient condition for the equivalence of the eigenpair between (\ref{eigp2}) and  (\ref{eigstt1}) is
\begin{equation}\label{condss1}
    \bm{C}^{\dag}\zeta={\bf{0}}\Longleftrightarrow\bm{Y}\bm{C}^{\dag}\zeta={\bf{0}}.
\end{equation}
The fundamental theory in linear algebra tells us that
(\ref{condss1}) holds if and only if
\begin{equation}\label{rankw}
    \textrm{rank}(\bm{C}^{\dag})=\textrm{rank}(\bm{Y}\bm{C}^{\dag}).
\end{equation}
By substituting (\ref{matrix12}) into (\ref{rankw}), we obtain the necessary and sufficient condition for the equivalence of the eigenpair between (\ref{eigp2}) and  (\ref{eigstt1}) is
\begin{equation}\label{auxiwww}
\textrm{rank}(\bm{M}^{\dag}\bm{Y}^{\dag})=\textrm{rank}(\bm{Y}\bm{M}^{\dag}\bm{Y}^{\dag}).
\end{equation}

We now prove that if $\d{\mu}^{\dag}=\d{\mu}>0$, then (\ref{auxiwww}) holds. In fact, it follows that $\bm{M}^{\dag}=\bm{M}>0$ from $\d{\mu}^{\dag}=\d{\mu}>0$. Hence, there exists an invertible square matrix $\bm{X}$ such that $\bm{M}=\bm{M}^{\dag}=\bm{X}^{\dag}\bm{X}$. In terms of the fundamental theory in linear algebra, it is known that
\begin{eqnarray}
&~~~&\textrm{rank}(\bm{M}^{\dag}\bm{Y}^{\dag})=\textrm{rank}(\bm{Y}^{\dag})=\textrm{rank}(\bm{X}\bm{Y}^{\dag})
\nonumber\\&&=\textrm{rank}(\bm{Y}\bm{X}^{\dag}\bm{X}\bm{Y}^{\dag})
    =\textrm{rank}(\bm{Y}\bm{M}^{\dag}\bm{Y}^{\dag}).\label{ranke1}
\end{eqnarray}
The rank equality (\ref{ranke1}) shows that the rank equality (\ref{auxiwww}) holds. As a result, we have already proved that the eigenpair between (\ref{eigp2}) and (\ref{eigstt1}) is equivalent provided that $\d{\mu}^{\dag}=\d{\mu}>0$.

According to the  above discussion, when the material is not magnetic lossy, i.e., under Case 1 and 2, each eigenpair between (\ref{eigp2}) and (\ref{eigstt1}) is equivalent. In this case, it is easy to show that each entry of eigenvector $\zeta$ in (\ref{eigstt1}) is the same.
Next, we briefly prove this conclusion.
In fact, since the matrix $\bm{M}$ is invertible, $\bm{Y}^{\dag}\zeta=\bf{0}$ can be derived from $\bm{C}^{\dag}\zeta=\bf{0}$ and $\bm{C}=\bm{Y}\bm{M}$. By means of (\ref{linearsol}), $\zeta=c\beta$ is valid, which shows that each entry of eigenvector $\zeta$ is the same.

For all the magnetic lossy materials, we cannot make sure each eigenvalue of (\ref{eigstt1}) is physical, because we cannot prove that (\ref{auxiwww}) is valid in this case. However, after we carry out many numerical experiments, these numerical results show that the augmented method is still free of all the spurious modes for certain magnetic lossy media.

\subsection{Projection Method}
In this subsection, we apply the projection method to solve resonant cavity problems under Case 1, 2, 3 and 4. Since the divergence-free condition is enforced in this numerical method, as a consequence, the projection method can remove all the spurious modes, including spurious zero modes.

It is known that all the solutions to (\ref{eigp2b}) form a linear subspace $\mathcal{V}$ in $\mathbb{C}^{n}$. Set
$\mathcal{V}=\textrm{span}\{q_{1},q_{2},\cdots,q_{r}\}$,
where $r=\dim{\mathcal{V}}$ and $q_{i}\in{\mathbb{C}^{n}}$. If the matrix $\bm{M}$ is invertible, then $\textrm{rank}(\bm{C})=\textrm{rank}(\bm{Y})=m-1$ is valid,
which implies $r=\dim{\mathcal{V}}=n-m+1$. Set $\bm{Q}=[q_{1},q_{2},\cdots,q_{r}]\in{\mathbb{C}^{n\times r}}$.

The basic idea of solving (\ref{eigp2}) is that choosing $\Lambda_{h}\in{\mathbb{C}}$ and $\xi\in{\mathcal{V}}$ such that
\begin{equation}\label{galerin}
    (\bm{A}\xi-\Lambda_{h}\bm{M}\xi)\perp \mathcal{V}.
\end{equation}
This is called the Galerkin condition \cite{Saad}. Set $\xi=\bm{Q}y$, where $y\in{\mathbb{C}^{r}}$. The Galerkin condition (\ref{galerin}) can be equivalently expressed the following equation
\begin{equation}\label{eigenp}
    (\bm{Q}^{\dag}\bm{A}\bm{Q})y=\Lambda_{h}(\bm{Q}^{\dag}\bm{M}\bm{Q})y.
\end{equation}
In order to compute the eigenpair of (\ref{eigenp}), the matrix $\bm{Q}$ must be given
in (\ref{eigenp}).

It is well-known that the popular numerical method of finding several eigenpairs of large-scale sparse matrices is an iterative method, since a fundamental operation in the iterative method is the matrix-vector multiplication and this operation is very efficient to the sparse matrix and vector. The efficient projection method needs to find a sparse matrix $\bm{Q}$ for the null space of the above-mentioned large sparse matrix $\bm{C}$. However, since Coleman and Pothen \cite{Coleman1986} prove that finding the sparsest basis for the null space of an underdetermined matrix is NP-complete hard, we {cannot} seek the sparsest matrix $\bm{Q}$ for the null space of $\bm{C}$.

Based on the importance of numerical stability, a set of normalized orthogonal bases $\{q_{1},q_{2},\cdots,q_{r}\}$ is used in this numerical calculation. In such a case,
$\bm{Q}^{\dag}\bm{Q}=\bm{I}_{r}$ holds, where $\bm{I}_{r}$ is the identity matrix of order $r$. Here, we employ singular value decomposition (SVD) technique to seek $\bm{Q}$. Suppose that
\begin{equation}\label{svd}
    \bm{C}=\bm{U}\bm{D}\bm{V}^{*}
\end{equation}
is the SVD of the matrix $\bm{C}$. We take $\bm{Q}=\bm{V}(:,n-r+1:n)$, where $\bm{V}(:,n-r+1:n)$ is a submatrix consisting of the last $r$ columns of $\bm{V}$, then $\bm{C}\bm{Q}=\bm{O}$ is valid.

If the eigenpair $(\Lambda_{h},y)$ of (\ref{eigenp}) is obtained by using the implicitly restarted Arnoldi methods \cite{Lehoucq}, then $(\Lambda_{h}, \bm{Q}y)$ is just an eigenpair of (\ref{eigp2}), which corresponds to a physical numerical mode of (\ref{eq:3}).

\subsection{The Advantage and Disadvantage among the Above Three Methods}
The main advantage of penalty method can preserve matrix size comparing with the augmented method. In addition, penalty method does not destroy the sparsity of the matrices comparing with the projection method. However, the main disadvantage is that penalty method will introduce spurious modes in solving 3-D closed cavity problem. In addition, selecting an appropriate penalty parameter $\alpha$ is not an easy work.

The main advantage of augmented method can preserve the sparsity of the matrix.
Moreover, when $\d{\mu}_{r}^{\dag}=\d{\mu}_{r}>0$, the augmented method is free of all the spurious modes. However, the main disadvantage of the augmented method is that the size of the matrix has increased.

The main advantage of projection method based on SVD can eliminate all the spurious modes even if the material is both electric and magnetic lossy.
However, the main disadvantage of projection method based on SVD is that this method is not efficient since the matrix $\bm{Q}$ is usually dense.

In a word, these three methods have their own advantages and disadvantages.

\section{Numerical Experiments}
In this section, we simulate three cavity problems by the above penalty method,
augmented method and projection method.
In order to distinguish the numerical eigenvalues associated with penalty method, augmented method and projection method, the numerical eigenvalues obtained by the penalty method, augmented method and projection method are denoted by $\Lambda_{h}(\textrm{pe},\alpha)$, $\Lambda_{h}(\textrm{au})$ and $\Lambda_{h}(\textrm{pr})$, respectively.
Here $\alpha$ is the parameter in the penalty method and it is usually a positive real number. It is worthwhile to point out that our adopted computational strategy is serial, instead of parallel.
\begin{table}[ht!]
\renewcommand{\arraystretch}{1.3}
\caption{\label{empty1} The Numerical Eigenvalues ($\Lambda_{h}$, $\mathrm{m}^{-2}$) Associated With the Dominant Mode from an Empty Spherical Resonant Cavity and CPU Time (Under Case 1)}
\centering
\begin{tabular}{cccccc}
 \hline
 $h(\textrm{m})$& 0.38493
 & 0.27062 & 0.22416 &0.16258 & Exact\\
 \hline
 $\Lambda_{h}(\textrm{pe},800)$&  7.71147 &7.62386 & 7.59006  &7.55655& 7.52793\\
 $t~(s)$& 10.3  & 18.5&  28.7 &40.6& --\\
 \hline
 $\Lambda_{h}(\textrm{au})$&7.71147 &7.62386 &  7.59006 &  7.55655&7.52793\\
  $t~(s)$& 12.7  & 21.8&  32.6 &60.7& --\\
 \hline
 $\Lambda_{h}(\textrm{pr})$& 7.71147  &7.62386  & 7.59006 &   7.55654&7.52793\\
  $t~(s)$&60.5   & 100.7 &  180.9 &360.8& --\\
 \hline
\end{tabular}
\end{table}
\subsection{Empty Spherical Resonant Cavity}

Let us consider an empty spherical resonator with the radius $r=1$\,m. The exact eigenvalue associated with the dominant mode is $\Lambda=7.52793\,\textrm{m}^{-2}$ \cite{jin2011theory}. Furthermore, the algebraic multiplicity of the exact eigenvalue $\Lambda$ is 3.
Suppose that the numerical eigenvalues $\Lambda_{h}^{(1)}$, $\Lambda_{h}^{(2)}$ and $\Lambda_{h}^{(3)}$ are the approximation of the exact eigenvalue $\Lambda$.
Set $\Lambda_{h}=(\Lambda_{h}^{(1)}+\Lambda_{h}^{(2)}+\Lambda_{h}^{(3)})/3$.
We employ the penalty method, augmented method and projection method to
solve this spherical resonant cavity problem, and then list the numerical
eigenvalues $\Lambda_{h}(\textrm{pe},\alpha)$, $\Lambda_{h}(\textrm{au})$ and $\Lambda_{h}(\textrm{pr})$ in Table \ref{empty1}. In order to compare with the efficiency of these three numerical methods, the CPU time is also given in Table \ref{empty1}.

\begin{table*}[ht]
\renewcommand{\arraystretch}{1.3}
\centering
\caption{\label{ser1} The Eigenvalues $\Lambda_{h}(\mathrm{au})$ and $\Lambda_{h}(\mathrm{pr})$ ($\mathrm{m}^{-2}$) with Physical Significance from Cylindrical Cavity (Under Case 2)}
\begin{tabular}{cccccc}
 \hline
 $h(\textrm{m})$& $0.1043$  & $0.0714$ & $0.0580$ &$0.0428$ &COMSOL\\
 \hline
 \hline
 $\Lambda_{h}^{1}(\textrm{au})$&  $24.0200 +11.9858\textrm{j}$ &$23.8807 +11.9245\textrm{j}$ & $23.8547 +11.9137\textrm{j}$  &$23.8226 +11.9097\textrm{j}$&$23.8230 +11.9085\textrm{j}$\\
  $\Lambda_{h}^{1}(\textrm{pr})$&  $24.0200+11.9858\textrm{j}$ &$23.8807 +11.9245\textrm{j}$ & $23.8547 +11.9137\textrm{j}$  &$23.8225 +11.9096\textrm{j}$&$23.8230 +11.9085\textrm{j}$\\
 \hline
 $\Lambda_{h}^{2}(\textrm{au})$& $26.6677 +13.3087\textrm{j}$&$26.4780 +13.2215\textrm{j}$ &  $26.4408 +13.2050\textrm{j}$ &  $26.3976+13.1853\textrm{j}$&$26.3968+13.1848\textrm{j}$\\
  $\Lambda_{h}^{2}(\textrm{pr})$& $26.6677+13.3087\textrm{j}$&$26.4780 +13.2215\textrm{j}$ &  $26.4408 +13.2050\textrm{j}$ &  $26.3974+13.1850\textrm{j}$&$26.3968+13.1848\textrm{j}$\\
  \hline
 $\Lambda_{h}^{3}(\textrm{au})$& $38.6158 + 0.0559\textrm{j}$  &$37.9265 + 0.0253\textrm{j}$  & $37.7824 + 0.0168\textrm{j}$ &   $37.6098+ 0.0079\textrm{j}$&$37.6067 + 0.0069\textrm{j}$\\
  $\Lambda_{h}^{3}(\textrm{pr})$& $38.6158 + 0.0559\textrm{j}$  &$37.9265 + 0.0253\textrm{j}$  & $37.7824 + 0.0168\textrm{j}$ &   $37.6097+ 0.0077\textrm{j}$&$37.6067 + 0.0069\textrm{j}$\\
\hline
\end{tabular}
\end{table*}

\begin{table*}[ht]
\renewcommand{\arraystretch}{1.3}
\centering
\caption{\label{ser2} The Eigenvalues $\Lambda_{h}(\mathrm{pe},1000)$, $\Lambda_{h}(\mathrm{au})$ and $\Lambda_{h}(\mathrm{pr})$ ($\mathrm{m}^{-2}$) with Physical Significance from Cylindrical Cavity (Under Case 4)}
\begin{tabular}{cccccc}
 \hline
 $h(\textrm{m})$& 0.1043  & 0.0714 & 0.0580 &0.0428 &COMSOL\\
 \hline
 \hline
  $\Lambda_{h}^{1}(\textrm{pe},1000)$&  $24.5131-7.5590\textrm{j}$ &$24.3324- 7.5554\textrm{j}$ & $24.2948-7.5597\textrm{j}$  &$24.2497 - 7.5591\textrm{j}$&$24.2476 - 7.5597\textrm{j}$\\
 $\Lambda_{h}^{1}(\textrm{au})$&  $24.5131 - 7.5590\textrm{j}$ &$24.3324 - 7.5554\textrm{j}$ & $24.2948 - 7.5597\textrm{j}$  &$24.2497-7.5591\textrm{j}$&$24.2476 - 7.5597\textrm{j}$\\
  $\Lambda_{h}^{1}(\textrm{pr})$&  $24.5131 - 7.5590\textrm{j}$ &$24.3324 - 7.5554\textrm{j}$ & $24.2947 - 7.5595\textrm{j}$  &$24.2495-7.5589\textrm{j}$&$24.2476 - 7.5597\textrm{j}$\\
 \hline
  $\Lambda_{h}^{2}(\textrm{pe},1000)$& $25.5404-9.7698\textrm{j}$&$25.3498- 9.7393\textrm{j}$ &  $25.3117- 9.7345\textrm{j}$ &  $25.2695 - 9.7261\textrm{j}$&$25.2649- 9.7244\textrm{j}$\\
 $\Lambda_{h}^{2}(\textrm{au})$& $25.5404 - 9.7698\textrm{j}$&$25.3498 - 9.7393\textrm{j}$ &  $25.3117 - 9.7345\textrm{j}$ &  $25.2695-9.7261\textrm{j}$&$25.2649- 9.7244\textrm{j}$\\
  $\Lambda_{h}^{2}(\textrm{pr})$& $25.5404- 9.7698\textrm{j}$&$25.3498-9.7393\textrm{j}$ &  $25.3123 - 9.7351\textrm{j}$ &  $25.2690 - 9.7259\textrm{j}$&$25.2649- 9.7244\textrm{j}$\\
  \hline
\end{tabular}
\end{table*}
 In these three numerical methods, the time and memory consumed by the projection method is the largest since the dense matrix $\bm{Q}$ obtained by SVD is used. This shows that the projection method is not efficient. In addition, the CPU time and memory consumed by the penalty method and the augmented method are roughly equivalent.
\begin{figure}[ht]
\centering
\includegraphics[width=8.8cm]{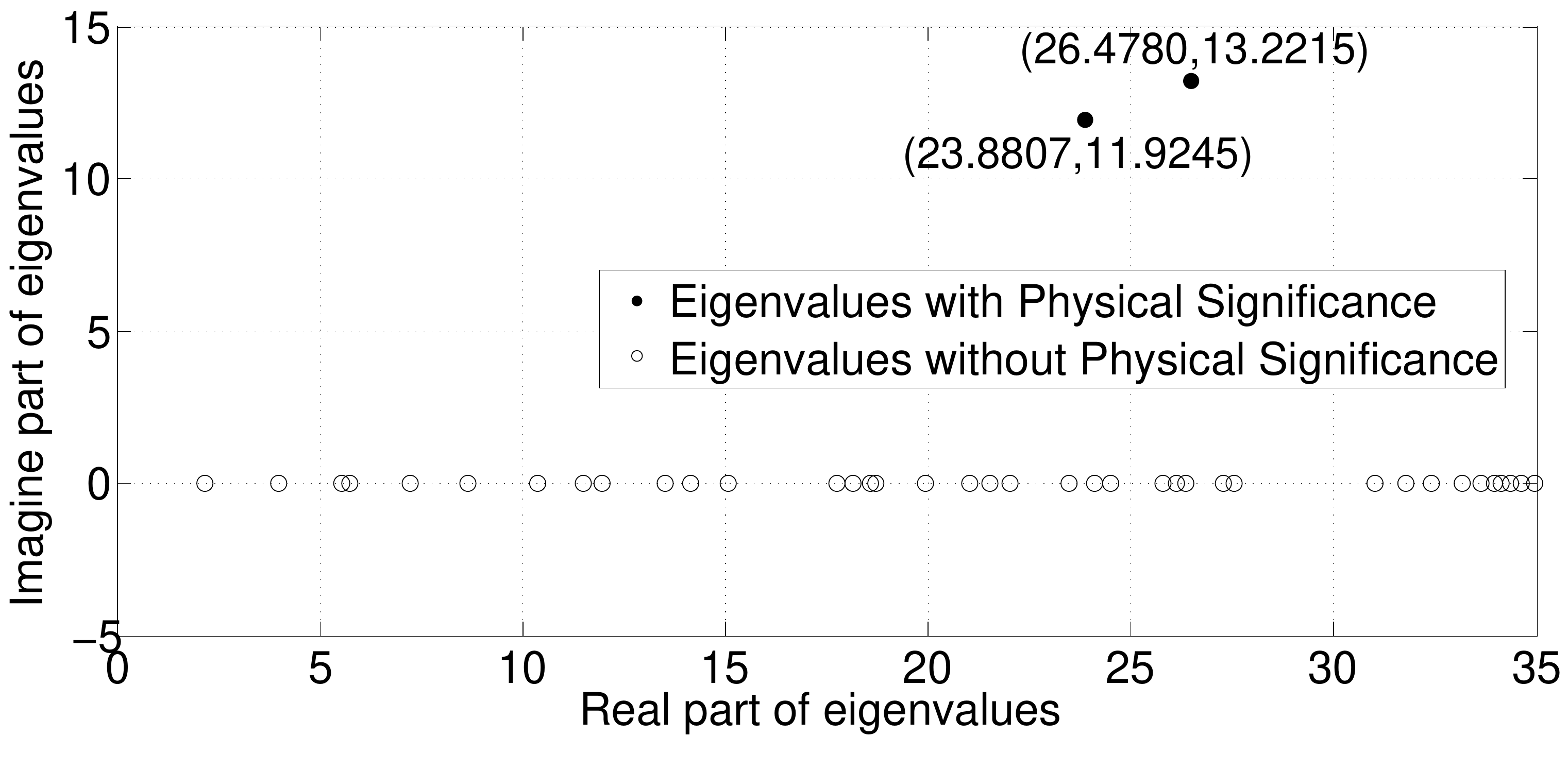}
\caption{Under the second mesh ($h=0.0714$\,m),
the eigenvalues associated with physical modes and spurious modes obtained by the penalty method with $\alpha=800$.}\label{fig2}
\end{figure}

In Table \ref{empty1}, we can also see that $\Lambda_{h}(\textrm{pe},800)=\Lambda_{h}(\textrm{au})\approx\Lambda_{h}(\textrm{pr})\approx\Lambda$ under the finest mesh. This shows that our numerical implementations are correct.
In this example, we do find that there are many eigenvalues associated with the spurious modes. These numerical eigenvalues are less than $\Lambda$ provided that the penalty parameter $\alpha$ is less than $700$. However, the numerical eigenvalues obtained by the augmented method and projection method are all physical.

\subsection{The Resonant Cavity Filled With Electric Lossy Media}

In this subsection we consider a cylindrical cavity with the radius $r=0.2$\,m and the height $h=0.5$\,m. Assuming that the relative permittivity and permeability tensor of the medium in the whole cylindrical cavity are
\begin{equation*}
    \d{\ep}_{r}=
    \begin{bmatrix}
    2-\textrm{j}&0&0\\
    0&2-\textrm{j}&0\\
    0&0&2
    \end{bmatrix},\quad
    \d{\mu}_{r}=
    \begin{bmatrix}
        2&-0.375\textrm{j}&0\\
    0.375\textrm{j}&2&0\\
    0&0&2
    \end{bmatrix}.
\end{equation*}

The first three numerical eigenvalues $\Lambda_{h}^{i}(\textrm{au})$ and $\Lambda_{h}^{i}(\textrm{pr})$ ($i=1,2,3$) are shown in Table \ref{ser1}. In Table \ref{ser1}, it can be observed that $\Lambda_{h}^{i}(\textrm{au})\approx\Lambda_{h}^{i}(\textrm{pr})$, $i=1,2,3$,
and they coincide with the eigenvalues corresponding to physical modes from COMSOL Multiphysics 5.2a. In the COMSOL simulation, the eigenvalues associated with physical modes are obtained by the fourth mesh ($h=0.0428$\,m). Notice that there are many spurious zero modes in the numerical results of COMSOL Multiphysics 5.2a. However, there are no any spurious modes in the numerical results of the augmented method and projection method.

Under the second mesh ($h=0.0714$\,m),
we employ the penalty method to solve this cylindrical cavity problem,
where $\alpha=800$ is taken in the numerical calculation, and then list the first forty numerical eigenvalues in Fig. \ref{fig2}. In Fig. \ref{fig2}, one can see that there are only two eigenvalues with physical significance, and the rest are all the eigenvalues without physical significance, whose imaginary part is zero.
Furthermore, these two eigenvalues with physical significance obtained by the penalty method
are equal to the ones obtained by augmented method.
In addition, in (\ref{eigstt1}),
we do find that the eigenvector $\zeta$ associated with each eigenvalue is a vector consisting of the same entry.

\subsection{The Resonant Cavity Filled With both Electric and Magnetic Lossy Media}
In this subsection we try to find the eigenmode of resonant cavity filled with an electric and magnetic lossy medium.
The penalty method, augmented method and projection method are used to solve this problem, and then we list the numerical eigenvalues associated with the first two physical modes in Table \ref{ser2}.

Suppose that the geometric shape of the cavity in this example is the same as the one in the example \emph{B}. In the cylindrical cavity, the relative permittivity and permeability tensor of the medium are
\begin{equation*}
    \d{\ep}_{r}=
    \begin{bmatrix}
    2+\textrm{j}&0&0\\
    0&2+\textrm{j}&0\\
    0&0&2
    \end{bmatrix},~~
    \d{\mu}_{r}=
    \begin{bmatrix}
        2-\textrm{j}&0.375\textrm{j}&0\\
    0.375\textrm{j}&2-\textrm{j}&0\\
    0&0&2
    \end{bmatrix}.
\end{equation*}
Obviously, the above material is both electric and magnetic lossy. Since the exact solution to this problem is unknown, we employ COMSOL Multiphysics 5.2a to simulate this problem, and then obtain the approximate eigenvalues of certain accuracy. The eigenvalues with physical significance from COMSOL are obtained by the fourth mesh ($h=0.0428$\,m). Notice that many spurious zero modes appear in the numerical results of COMSOL.


The numerical eigenvalues from the penalty method and projection method are listed in Table \ref{ser2}. In Table \ref{ser2}, we can see that $\Lambda_{h}^{i}(\textrm{pe},1000)\approx\Lambda_{h}^{i}(\textrm{au})\approx\Lambda_{h}^{i}(\textrm{pr})$, $i=1,2$, which coincide with the eigenvalues corresponding to physical significance from COMSOL. Here it is worthwhile to emphasize that the projection method can remove all the spurious modes.


\section{Conclusion}
The finite element method can be applied to solve 3-D closed cavity problem filled with anisotropic and nonconductive media. The matrix system resulting from the finite element method is a constrained generalized eigenvalue problem. This difficult problem can be solved by the penalty method, augmented method and projection method. The penalty method {cannot} remove all the spurious modes. We prove that the augmented method is free of all the spurious modes if the medium is not magnetic lossy.
When the medium is both electric and magnetic lossy, the projection method based on SVD technique can deal with this type of resonant cavity problem very well. However, the projection method based on SVD technique is not efficient. In future, we would like to give an efficient iterative method to solve the constrained generalized eigenvalue problem.
\section*{Acknowledgement}
We gratefully acknowledge the help of Prof. Qing Huo Liu, who has offered us valuable suggestions in the revision and Dr. Yuanguo Zhou for improving our English writing.


\begin{IEEEbiography}[{\includegraphics[width=1.0in,height=1.25in,clip]{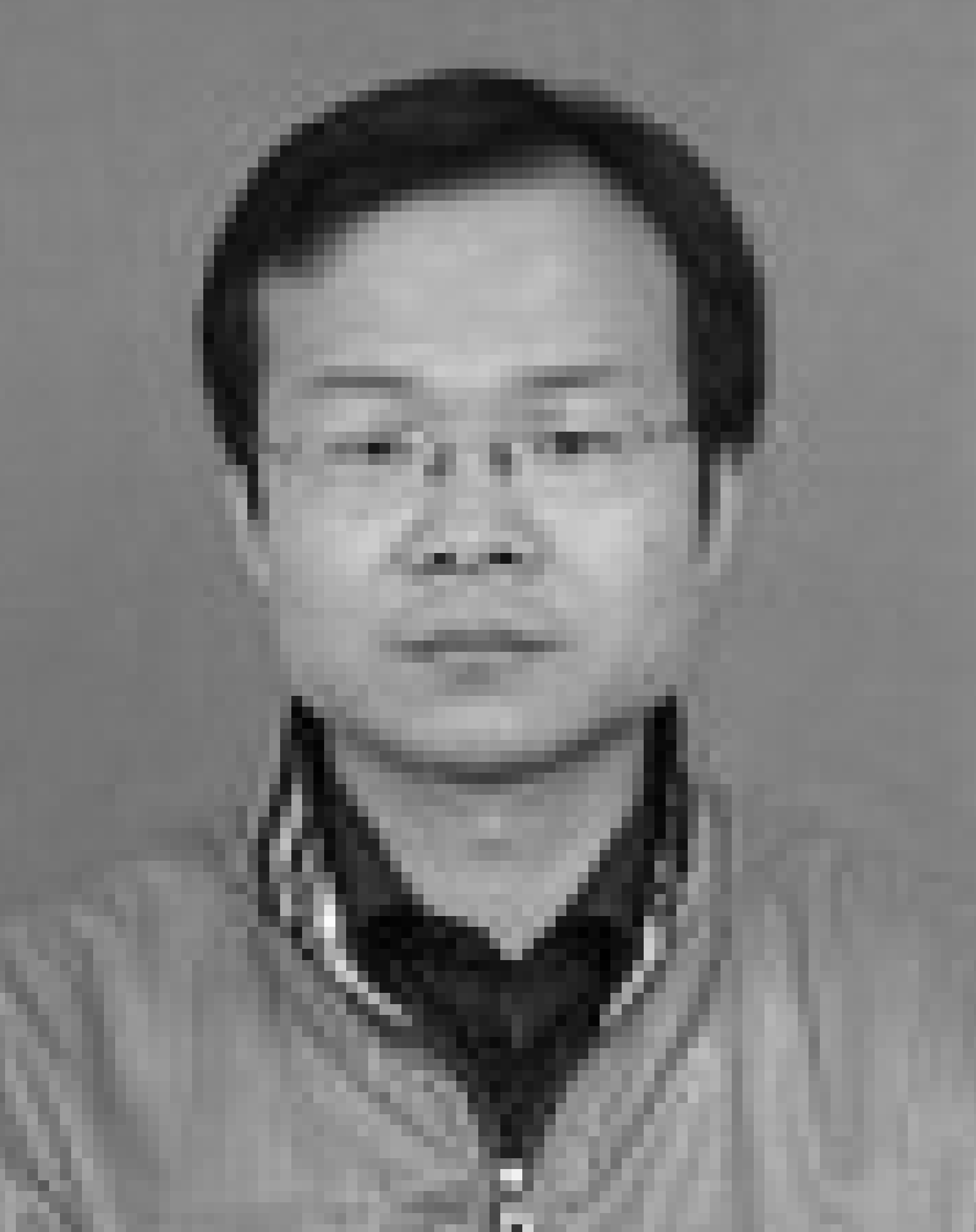}}]{Wei Jiang} was born in Wuhan, Hubei province, China, in 1984. He received the B.S. degree in mathematics and applied mathematics from Hubei Normal University, Huangshi, China, in 2008, the M.S. degree in computational mathematics from Guizhou Normal University, Guiyang, China, in 2011, and the Ph.D. degree in radio physics from Xiamen University, Xiamen, China, in 2016.

From September 2016 to July 2018, he has engaged in postdoctoral research at the Institute of Geophysics and Geomatics, China University of Geosciences, Wuhan, China. He is currently a lecturer at the School of Mechatronics Engineering, Guizhou Minzu University, Guiyang, China. Until now, he has published over 10 papers in refereed journals. His current research interests include the finite-element method for partial differential equations, applied numerical algebra, and computational electromagnetics, especially for eigenvalue problems in electromagnetics.
\end{IEEEbiography}

\begin{IEEEbiography}[{\includegraphics[width=1.0in,height=1.25in,clip]{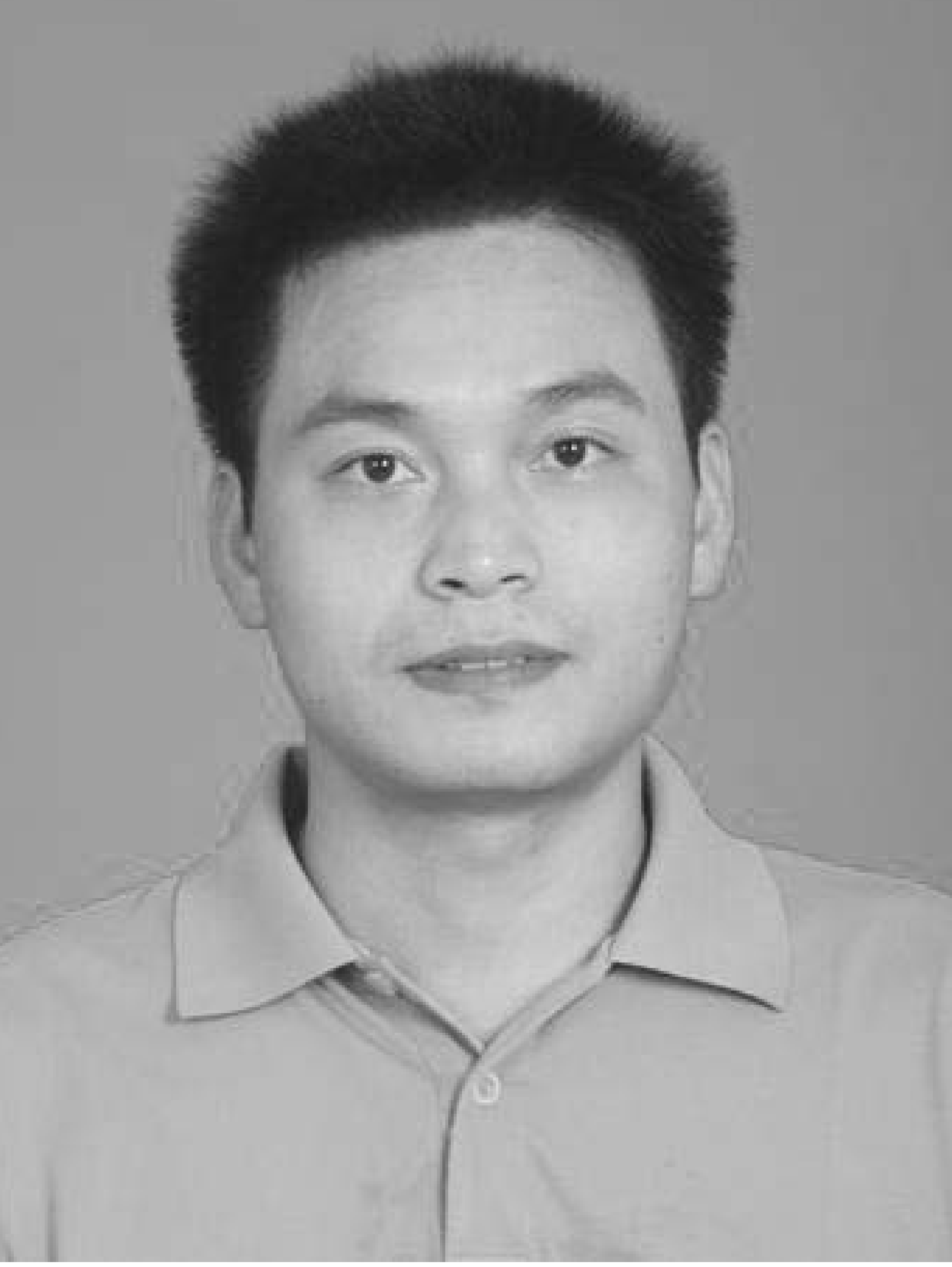}}]{Jie Liu} was born in Tongren, Guizhou province, China, in 1984. He received the B.S. degree in mathematics
and applied mathematics and the M.S. degree in computational mathematics from Guizhou Normal University, Guiyang, China, in 2008 and 2012, respectively,
and the Ph.D. degree in electromagnetic
fields and microwave techniques from the Institute of
Electromagnetics and Acoustics, Xiamen University,
Xiamen, China, in 2019.

From July 2012 to July 2015, he was a Lecturer
with the School of Mathematics and Statistics,
Guizhou University of Finance and Economics,
Guiyang. He is currently a Post-Doctoral Researcher with the School of
Informatics, Xiamen University, Xiamen, China. His current research interests
include the finite-element method and the spectral-element method for partial
differential equation eigenvalue problems and computational electromagnetics.
\end{IEEEbiography}
\end{document}